\input amstex
\mag=\magstep1
\documentstyle{amsppt}

\define\cc{{\Bbb C}}
\define\ccnz{{\Bbb C^n - \{0\}}}
\define\cz{{\Bbb C^\times}}

\define\dz{\frac{d}{dz}}
\define\dzn{{dz_1}^\nu\cdots{dz_n}^\nu}
\define\ep{\varepsilon}
\define\ev{\operatorname{ev}}
\define\ez{{e_0}}
\define\fai{{\Cal F}}
\define\faig{{\Cal F^G}}
\define\gfg{{\frak g}}
\define\Hom{\operatorname{Hom}}
\define\lz{{L_0}}

\define\om{\Cal O_M}
\define\ome{\Cal O_M(E)}
\define\pr{\operatorname{pr}}
\define\sgns{(\operatorname{sign}\sigma)}
\define\sn{{\frak S_n}}
\define\tconv{{\text{conv.}}}
\define\talg{{\text{alg.}}}
\define\talt{{\text{alt.}}}
\define\tsym{{\text{sym.}}}
\define\ww{{W_1}}
\define\wnzo{1_{\nu_0}\otimes}
\define\wnzotln{1_{\nu_0}\otimes\Cal O_{\Bbb C^n}(\tau^{\lambda_1,\dots,\lambda_n}_{\nu_1,\dots,\nu_n})}
\define\wnzotlnn{1_{\nu_0}\otimes\Cal O_{\Bbb C^n}((\tau^{\lambda}_{\nu})^n)}
\define\wnzoztn{1_{\nu_0}\otimes\Cal O_{(\Bbb C^\times)^n}(\tau_{\nu_1,\dots,\nu_n})}
\define\wnztlnn{1_{\nu_0}\otimes(T^{\lambda}_{\nu})^{\otimes n}}
\define\wnztnztln{1_{\nu_0}\otimes\bigotimes^n_{i=1}(T^\times_{\nu_i}/T^{\lambda_i}_{\nu_i})}
\define\wnztnztlnn{1_{\nu_0}\otimes(T^\times_{\nu}/T^{\lambda}_{\nu})^{\otimes n}}

\topmatter
\title An application of the second Riemann continuation theorem\\
 to cohomology of the Lie algebra of vector fields\\
 on the complex line. \endtitle
\date{August 21, 2005}\enddate
\rightheadtext{Cohomology of the Lie algebra of vector fields 
 on the line}
\author Nariya KAWAZUMI
\endauthor
\affil
Department of Mathematical Sciences,\\ University of Tokyo
\endaffil
\address
Tokyo, 153-8914 Japan 
\endaddress
\email
kawazumi\@ms.u-tokyo.ac.jp
\endemail

\abstract{
We study cohomology groups of the Lie algebra of vector fields 
on the complex line, $W_1$, with values in the tensor fields 
in several variables.  From a generalization by Scheja 
of the second Riemann (Hartogs) continuation theorem, 
we deduce a cohomology exact sequence of the subalgebra of 
$W_1$ consisting of vectors having a zero at the origin.  
As applications, we compute the cohomology algebra of $W_1$ with values 
in the functions on $\Bbb C^n$ explicitly, and establish 
a certain vanishing theorem for the cohomology of $W_1$ 
with values in the quadratic differentials 
in several variables, which is closely related 
to the moduli space of Riemann surfaces.}
\endabstract
\endtopmatter
\document
\heading Introduction. \endheading
The topological Lie algebra of complex analytic vector fields 
on an open Riemann surface $O$, $L(O)$, acts continuously on 
the space of complex analytic tensor fields on the product 
space $O^n$ by the diagonal Lie derivative action. 
By a complex analytic analogue [Ka1] 
of the Bott-Segal addition theorem [BS], 
the computation of the (continuous) cohomology group of $L(O)$ 
with values in the space of tensor fields on $O^n$ 
was reduced to that for the case 
where $O$ is the complex line $\Bbb C$ 
together with a topological study of certain sheves on the space $O^n$. 
The purpose of the present paper is to provide a geometric tool 
for the computation for this case $O = \Bbb C$. \par
We fix our notations. 
Let $\ww := L(\cc) = H^0(\cc; \Cal O_\cc(T\cc))$ be 
the Lie algebra of complex analytic vector fields 
on the complex line $\cc$ with the topology of uniform convergence 
on compact sets.  
The closed subalgebra $\lz := \{ X \in \ww; X(0) = 0 \}$ of $\ww$ 
often plays more important roles than $\ww$ itself. 
We recall two kinds of $\lz$ modules.
\roster
\item Let $\nu \in \Bbb Z$. 
The Lie algebra $\lz$ acts on the $1$ dimensional complex vector space 
$1_\nu = \cc1_\nu$ with the preferred base $1_\nu$ by 
$$
(\xi(z)\dz)\cdot 1_\nu = \nu\xi'(0)1_\nu \quad (\xi(z)\dz \in \lz).
$$
The $\lz$ module $1_\nu$ is naturally isomorphic to the $\nu$-cotangent 
space $(T^*_0\cc)^{\otimes\nu}$ at the origin.
\item For $\nu$ and $\lambda \in \Bbb Z$ 
we denote by $T^\lambda_\nu$ the Fr\'echet space of 
the meromorphic $\nu$-covariant tensor fields on $\cc$ 
with a pole only at the origin of order $\leq \lambda$, 
$T^\lambda_\nu = H^0(\cc; \Cal O_\cc((T^*\cc)^{\otimes\nu}\otimes 
[0]^{\otimes \lambda}))$. 
Here $[0]$ is the line bundle induced by the divisor $0 \in \cc$. 
The algebra $\lz$ acts on $T^\lambda_\nu$ by the Lie derivative
$$
(\xi(z)\dz)\cdot(f(z)dz^\nu) = (\xi(z)f'(z) + \nu\xi'(z)f(z))dz^\nu,
$$
where $\xi(z)\dz \in \lz$ and $f(z)dz^\nu \in T^\lambda_\nu$.\par
When $\lambda = 0$, $T^0_\nu$ is a $\ww$ module, 
which is denoted by $T_\nu$ for simplicity. 
The $\ww$ module of $\nu$-covariant tensor fields on $\cz = \cc - \{0\}$ 
is denoted by $T^\times_\nu$, i.e., 
$T^\times_\nu := H^0(\cc - \{0\}; \Cal O_\cc((T^*\cc)^{\otimes\nu}))$.
\endroster
Our purpose is to calculate the cohomology groups
$$
H^*(\ww; {\bigotimes}^n_{i=1}T_{\nu_i}) \quad\text{and}\quad 
H^*(\lz; 1_{\nu_0}\otimes{\bigotimes}^n_{i=1}T^{\lambda_i}_{\nu_i}),
$$
for arbitrary integers $\nu_0, \nu_1, \dots. \nu_n$ 
and $\lambda_1, \dots, \lambda_n$. 
Here and throughout this paper $\otimes$ means the completed tensor 
product over the complex numbers $\cc$. 
The Shapiro Lemma (Remark 2.4) implies the isomorphism 
$H^*(\ww; \bigotimes^n_{i=1}{T_{\nu_i}}) \cong 
H^*(\lz; 1_{\nu_1}\otimes\bigotimes^n_{i=2}{T_{\nu_i}})$, 
so that the calculation of the former is reduced to that of the latter.\par
First of all the cohomology group $H^*(\ww; T_\nu) \cong H^*(\lz; 1_\nu)$ 
was determined by Goncharova [Go]. 
Studying a filtration of the space 
$1_{\nu_0}\otimes T^{\lambda_1}_{\nu_1}$ derived from 
the action of the vector $\ez := z\dz \in \lz$ in detail, 
Feigin and Fuks [FF] determined the cohomology group 
$H^*(\lz; 1_{\nu_0}\otimes T^{\lambda_1}_{\nu_1})$ completely. 
In the present paper we study the cohomology group 
$H^*(\lz; 1_{\nu_0}\otimes\bigotimes^n_{i=1}T^{\lambda_i}_{\nu_i})$ for the case $n \geq 2$ in a geometric way. 
Using a generalization by Scheja [S] of the second Riemann (Hartogs) 
continuation theorem, we introduce a cohomology exact sequence (3.2)
$$
\multline
\cdots\to H^{q-n}(\lz; \wnztnztln) \overset{d_n}\to\to 
H^q(\lz; 1_{\nu_0}\otimes\bigotimes^n_{i=1}T^{\lambda_i}_{\nu_i})\\
\to H_\lz^q(\ccnz; \wnzotln)\to\cdots,
\endmultline
$$
which we call {\it the fundamental exact sequence for the $\lz$ module 
$1_{\nu_0}\otimes\bigotimes^n_{i=1}T^{\lambda_i}_{\nu_i}$}. 
The $\ez$-invariants $C^*(\lz; \wnztnztln)^\ez$ is of finite dimension,
and so is the first term.  
The third term is the $\lz$ equivariant cohomology of the space $\ccnz$ 
with values in the sheaf of $\lz$ modules $\wnzotln$ (\S2). 
This term is reduced to the cohomology of $\lz$ with values 
in the tensor fields in fewer variables by the main result 
of our previous paper [Ka1] Theorem 5.3.  
In \S1 we give an exposition of the equivariant cohomology theory 
for a Lie algebra.\par
Three results are deduced from the sequence (3.2).
\roster
\item The cohomology of $\lz^\talg := z\cc[z]\dz$ with values 
in the space of {\it algebraic} tensor fields on the complex torus 
$(\cc - \{0\})^n$ coincides with that of $\lz$ with values 
in the space of {\it complex analytic} tensor fields of the same type 
on the torus (Corollary 3.5). 
\item The cohomology of $\ww$ with values in the algebra 
of complex analytic functions on $\cc^n$, or equivalently, 
that of $\ww^\talg = \cc[z]\dz$ with values in the algebra of 
polynomials in $n$ variables, $\cc[z_1, \dots, z_n]$, 
is determined explicitly (Corollary 3.7).
\item A vanishing theorem (Theorem 4.7) on the quadratic differentials 
in $p$ variables:
$$
H^q(\ww; {\bigwedge}^pQ) = 0, \quad \text{if } p > q,
$$
where $Q := T_2$ and $\bigwedge^p$ denotes the completed 
$p$-fold alternating tensor product on the complex numbers $\cc$.  
This suggests that the Virasoro equivariant $(p, q)$ cohomology of 
the dressed moduli of compact Riemann surfaces would vanish 
for $p > q$ [Ka2].
\endroster
\par
The author would like to express his own gratitude to Yukio Matsumoto 
for constant encouragement and to Shigeyuki Morita, Kazushi Ahara and 
Masanori Kobayashi for helpful discussions.\par
This paper is a revised version of an unpublished preprint 
written in 1993 (University of Tokyo, UTMS 93-18). \par
\demo{Contents}\par
\S1. Equivariant cohomology.\par
\S2. Preliminaries.\par
\S3. Fundamental exact sequence.\par
\S4. Action of the symmetric group.\par
\enddemo

\beginsection 1. Equivariant cohomology.\par
We begin by an exposition on the equivariant cohomology 
for a Lie algebra. 
Let $\gfg$ be a complex topological Lie algebra.  
A $\gfg$ {\it module} means a complex topological vector space 
on which $\gfg$ acts continuously.  
The standard continuous cochain complex of 
the topological Lie algebra $\frak g$ 
with values in a $\frak g$ module $N$ is denoted by
$$
C^*(\frak g; N) = {\bigoplus}_{p\ge0}C^p(\frak g; N).
$$
Here $C^p(\frak g; N)$ is the linear space of continuous 
alternating multi-linear mappings $c: 
\undersetbrace{p}\to{\gfg\times\cdots\times\gfg} \to N$. 
They can be regarded as continuous linear maps of the completed 
alternating tensor product $\bigwedge^p\gfg$ into $N$. 
Hence we may identify $C^p(\gfg; N)$ with $\Hom(\bigwedge^p\gfg, N)$. 
The cohomology group of the complex $C^*(\frak g; N)$ is 
called the (continuous) cohomology group of $\frak g$ with values 
in $N$ and denoted by $H^*(\frak g; N)$.  
When $N$ is the trivial $\frak g$ module $\Bbb C$, 
we abbreviate them to $C^*(\frak g)$ and $H^*(\frak g)$ respectively. 
For details, see, for example, [HS].\par
If $\fai$ is a sheaf of $\gfg$ modules over a topological space $M$, 
the cochain complex of sheaves over $M$
$$
C^*(\gfg; \fai): U \overset{\text{open}}\to{\subset} M 
\mapsto C^*(\gfg; \fai(U))
$$
is defined. 
We denote by $H^*_{\gfg}(M; \fai)$ the hypercohomology group 
of the cochain complex of sheaves over $M$ with respect 
to the functor $\Gamma(M; \cdot)$ 
(= the sections of the sheaf $\cdot$ over $M$) ([G] ch.0, \S11.4, pp.32-) 
and call it {\it the $\gfg$ equivariant cohomology group of $M$ 
with values in the sheaf $\fai$}. Equivalently we define 
$$
H^*_{\gfg}(M; \fai) := H^*(\operatorname{Total}(\Gamma(M; \Cal C^{*, *})))
$$
for an injective right Cartan-Eilenberg resolution 
$\Cal C^{*, *} = (\Cal C^{i, j})_{i, j \geq 0}$ 
of the complex $C^*(\gfg; \fai)$ (cf. ibid\.loc\.cit\.). 
Needless to say this cohomology is 
an analogue of the equivariant cohomology 
of a space $M$ acted on by a transformation group $G$, i.e., 
the ordinary cohomology of the homotopy quotient $M_G := E_G\times_GM$.\par
There exist two spectral sequences converging 
to the equivariant cohomology $H^*_{\gfg}(M; \fai)$
$$
\split
`E^{p, q}_2 & = H^p(H^q(M; C^*(\gfg; \fai)))\\
``E^{p, q}_2 & = H^p(M; H^q(\gfg; \fai)),
\endsplit\tag1.1
$$
where we denote by $H^*(\gfg; \fai)$ the sheaf over $M$ 
given by the presheaf
$$
U \overset{\text{open}}\to{\subset} M \mapsto H^*(\gfg; \fai(U)).
$$
\par
To study the $`E_2$ term, we introduce a condition (1.2) related 
to a locally convex space $F$:\roster
\item $F$ is a Fr\'echet nuclear space,
\item there exists a projective system 
$\{ F_i, u^{i+1}_i: F_{i+1} \to F_i\}^\infty_{i=1}$ of 
the {\it strong duals} $F_i$ of Fr\'echet nuclear spaces 
and compact maps $u^{i+1}_i: F_{i+1} \to F_i$ 
satisfying a topological isomorphism $F = \varprojlim F_i$, and
\item the natural projection $F \to F_i$ has a dense image in each $F_i$,
\endroster
(cf\. [K]). Then we have ([Sw])
\proclaim{Lemma 1.3} {\rm (K\"unneth formula)} 
Let $F$ be a locally convex space satisfying the condition (1.2) 
and $A^*$ a Fr\'echet cochain complex with 
$\dim H^q(A^*) < \infty$ for each $q$. 
Then we have 
$$
H^*(\Hom(F, A^*)) = \Hom(F, H^*(A^*)),
$$
where $\Hom$ means the continuous linear maps.\endproclaim
Now we confine ourselves to the case $M$ is a {\it finite dimensional} 
complex analytic manifold and $\fai = \ome =$ the sheaf of germs 
of sections of a complex analytic vector bundle $E$ over $M$. 
Let $F$ be a locally convex space satisfying the condition (1.2). 
We denote by $A^q(U; E)$ the Fr\'echet space of $C^\infty$ sections 
of $E\otimes\bigwedge^q\overline{T^*M}$ on an open subset $U$ of $M$. 
Consider the sheaf $\Hom(F, \ome)$ (resp\. $\Hom(F, \Cal A^q(E))$) 
given by
$$
\split
& U \overset{\text{open}}\to{\subset} M \mapsto \Hom(F, \ome(U)),\\
(\text{resp. } \quad & U \overset{\text{open}}\to{\subset} M 
\mapsto \Hom(F, A^q(U; E)),)
\endsplit
$$
where we endow $\ome(U)$ (resp\. $A^q(U; E)$) 
with the Fr\'echet topology of uniformly convergence 
on compact subsets (resp\. with the $C^\infty$ topology).  
If $U \subset M$ is a Stein open subset, we have
$$
H^*(\Hom(F, A^*(U; E))) = \Hom(F, \ome(U)) \quad\text{(in dim. $0$)},\tag1.4
$$
by Lemma 1.3 applied to the acyclic augmented cochain complex 
$A^*(U;\allowmathbreak E)$ with the augmentation 
$\ome(U) \hookrightarrow A^0(U; E)$. \par
Passing to the inductive limit, we obtain a fine resolution of 
the sheaf $\Hom(F, \allowmathbreak \ome)$
$$
0 \to \Hom(F, \ome) \to \Hom(F, \Cal A^0(E)) \to \Hom(F, \Cal A^1(E)) 
\to \cdots.
$$
From (1.4) any Stein covering of $M$ is acyclic for the sheaf 
$\Hom(F, \ome)$. 

\proclaim{Proposition 1.5} If $\gfg$ satisfies the condition (1.2), 
and $M$, $E$ and $H^q(M; \ome)$ for $q \neq 0$ are all finite dimensional, 
then the natural map
$$
`E^{p, q}_2 = H^p(H^q(M; C^*(\gfg; \ome))) \to H^p(\gfg; H^q(M; \ome))
$$
is an isomorphism.
\endproclaim
\demo{Proof} Using Lemma 1.3 again, we have 
$$
H^*(M; \Hom(F, \ome)) = \Hom(F, H^*(M; \ome)).
$$
Substituting $F = \bigwedge^q\gfg$ to this isomorphism, we obtain
$$
H^q(M; C^p(\gfg; \ome)) = C^p(\gfg; H^q(M; \ome)),
$$
as was to be shown.\qed\enddemo
\demo{Example 1.6} Let $M$ be a compact K\"ahler manifold and 
$E$ the $n$-cotangent bundle $\bigwedge^nT^*M$. 
Since $\gfg$ acts on $H^*(M; \om(\bigwedge^nT^*M)) = H^{n, *}(M) 
\subset H^{n+*}(M)$ trivially, 
we have $`E^{p, q}_2 = H^p(\gfg)\otimes H^{n, q}(M)$. \enddemo
\demo{Example 1.7} Let $O$ be an open Riemann surface and $S$ 
a finite subset of $O$. We denote by $L(O,S)$ 
the Lie algebra of complex analytic vector fields on $O$ 
which have zeroes at all points in $S$. 
$L(O,S)$ is a Fr\'echet space satisfying the condition (1.2) 
with respect to the topology of uniform convergence on compact sets. 
Let $E \to M$ be a complex analytic vector bundle over 
a finite dimensional Stein manifold $M$. 
If the Lie algebra $L(O,S)$ acts on the sheaf of topological 
linear spaces $\om(E)$ continuously, then we have 
$$
H^*_{L(O,S)}(M; \om(E)) = H^*(L(O,S); \om(E)(M)).
$$
Hence we obtain a spectral sequence 
$$
``E^{p, q}_2 = H^p(M; \Cal H^q) \Rightarrow
H^{p+q}(L(O,S); \om(E)(M)),
$$ where $\Cal H^q$ is a sheaf over $M$ whose stalk at $x \in M$ is given by
$$
\Cal H^q_x = H^q(L(O,S); \om(E)_x),
$$
which we call the Re\v setnikov spectral sequence. 
See [R] and [Ka1] \S9.\enddemo
Let $E$ be a complex analytic vector bundle of finite rank over $\cc^n$. 
Consider the alternating \v Cech complex $C^*(\frak U) = 
C^*(\frak U; \Cal O_{\cc^n}(E))$ with respect to the $(n-1)$ dimensional 
Stein covering $\frak U = \{U_j\}^n_{j=1}$ of $\ccnz$ given by 
$$
U_j = \{(z_1, \dots, z_n) \in \Bbb C^n; z_j \neq 0\}. \tag1.8
$$
In view of a theorem of Scheja [S], 
$$
H^q(\ccnz; \Cal O_{\Bbb C^n}(E)) = H^q(\cc^n; \Cal O_{\cc^n}(E))
$$
for $q \leq n-2$.  
Applying Lemma 1.3 to the acyclic augmented complex $C^*(\frak U)$ 
with $C^{-1}(\frak U) = \ome(\cc^n)$ and 
$C^n(\frak U) = H^{n-1}(\ccnz; \Cal O_{\cc^n}(E))$, we have 
$$
H^*(\ccnz; \Hom(F, \Cal O_{\cc^n}(E))) = \Hom(F, H^*(\ccnz; \Cal O_{\cc^n}(E))).
$$
Consequently we obtain
\proclaim{Proposition 1.9} If $\gfg$ satisfies the condition (1.2), then the $`E_2$ term (1.1) converging to the equivariant cohomology $H^*_\gfg(\ccnz; \Cal O_{\cc^n}(E))$ is given by
$$
`E^{p, q}_2 =
\cases
H^p(\gfg; \Cal O_{\cc^n}(E)(\cc^n)), & \quad \text{ if } q = 0,\\
H^p(\gfg; H^{n-1}(\ccnz; \Cal O_{\cc^n}(E))), & \quad \text{ if } q = n - 1,\\
0, & \quad \text{ otherwise. }
\endcases
$$
\endproclaim

\beginsection 2. Preliminaries.\par
The topological Lie algebras $\ww$ and $\lz$, 
the $\lz$ modules $1_\nu$ and $T^\lambda_\nu$, 
and the $\ww$ modules $T_\nu$ and $T^\times_\nu$ 
are defined in Introduction. 
For $\lambda$ and $\nu \in \Bbb Z$ we denote by $\tau^\lambda_\nu$ 
the line bundle $(T^*\cc)^{\otimes\nu}\otimes[0]^{\otimes\lambda}$ 
over the complex line $\cc$, 
where $[0]$ is the line bundle induced by the divisor $0 \in \cc$.  
If $\lambda = 0$, we denote $\tau_\nu := \tau^0_\nu$.  
We have $T^\lambda_\nu = H^0(\cc; \Cal O_\cc(\tau^\lambda_\nu))$, 
$T_\nu = H^0(\cc; \Cal O_\cc(\tau_\nu))$ and 
$T^\times_\nu = H^0(\cc - \{0\}; \Cal O_\cc(\tau_\nu))$. 
We denote by $\tau^{\lambda_1,\dots,\lambda_n}_{\nu_1,\dots,\nu_n}$
the line bundle over $\cc^n$ given by 
$$
(\pr^*_1\tau^{\lambda_1}_{\nu_1})\otimes \cdots 
\otimes(\pr^*_n\tau^{\lambda_n}_{\nu_n}), \quad 
(\lambda_1, \dots, \lambda_n, \nu_1, \dots, \nu_n \in \Bbb Z),
$$
where $\pr_i: \cc^n \to \cc$ is the $i$-th projection. 
By the diagonal Lie derivative action, $\wnzotln$ 
($\nu_0 \in \Bbb Z$) 
is a sheaf of $\lz$ modules.  
From the nuclear theorem we have $\lz$ isomorphisms
$$
1_{\nu_0}\otimes{\bigotimes}^n_{i=1}T^{\lambda_i}_{\nu_i} 
= H^0(\cc^n; \wnzotln).
$$
If $\lambda_1 = \cdots = \lambda_n = \lambda$ and 
$\nu_1 = \cdots = \nu_n = \nu$, we abbreviate 
$(\tau^\lambda_\nu)^n := 
\tau^{\lambda_1,\dots,\lambda_n}_{\nu_1,\dots,\nu_n}$. 
\par
We recall a basic fact on an $\ez := z\dz$ invariant complex.
\proclaim{Lemma 2.1} For the $\lz$ modules 
$N = 1_{\nu_0}\otimes\bigotimes^n_{i=1}T^{\lambda_i}_{\nu_i}$ and 
$1_{\nu_0}\otimes\bigotimes^n_{i=1}
(T^\times_{\nu_i}/T^{\lambda_i}_{\nu_i})$, 
the inclusion $C^*(\lz; N)^\ez \subset C^*(\lz; N)$ 
induces a cohomology equivalence.
\endproclaim
\demo{Proof}(cf\. [Ka1]\S2). 
The action of the multiplicative group $\cz := \cc - \{0\}$ 
on the complex line
$$
T_t: \cc \to \cc, \quad z \mapsto tz, \quad (t \in \cz)
$$
induces the actions $T_t$ on $N$ and $\lz$ itself such that
$$
t\frac{d}{dt}T_t = T_t\ez = \ez T_t.\tag2.2
$$
Using the averaging operator 
$$
C^*(\lz; N) \to C^*(\lz; N)^{e_0}, \quad 
c\mapsto \int^{1}_{0}(T_{\exp2\pi\sqrt{-1}\theta}c)d\theta,
$$
we obtain the desired cohomology equivalence.\qed\enddemo
Now we introduce some variants of $\lz$, $\ww$ and $T^\lambda_\nu$. We define $\ww^\talg$, $\ww^\tconv$, $\lz^\talg$ and $\lz^\tconv$ by
$$
\split
& \ww^\talg := \cc[z]\dz, \quad\text{(the polynomial vector fields),}\\
& \ww^\tconv := \Cal O_\cc(T^*\cc)_0 = \cc\{z\}\dz, 
\quad\text{(the germs at $0 \in \cc$),}\\
& \lz^\talg := \{ X \in \ww^\talg; X(0) = 0 \},\quad\text{and}\\
& \lz^\tconv := \{ X \in \ww^\tconv; X(0) = 0 \}.
\endsplit
$$
Clearly $\ww^\talg \subset \ww \subset \ww^\tconv$ and 
$\lz^\talg \subset \lz \subset \lz^\tconv$.  ${T^\lambda_\nu}^\talg$ 
(resp\. ${T^\lambda_\nu}^\tconv$) defined by 
$$
\split
{T^\lambda_\nu}^\talg := \frac1{z^\lambda}\cc[z]dz^\nu \quad
\left(\text{resp.}\,\,\, {T^\lambda_\nu}^\tconv := 
\Cal O_\cc(\tau^\lambda_\nu)_0 = \frac1{z^\lambda}\cc\{z\}dz^\nu\right)
\endsplit
$$
is a $\lz^\talg$ (resp\. $\lz^\tconv$) module.  
If $\lambda = 0$, we denote ${T_\nu}^\talg := {T^0_\nu}^\talg$ and 
${T_\nu}^\tconv := {T^0_\nu}^\tconv$ for simplicity. 
We have no essential difference between these variants.
\proclaim{Lemma 2.3} We have natural isomorphisms
$$\align
& H^*(\lz; 1_{\nu_0}\otimes{\bigotimes}^n_{i=1}{T^{\lambda_i}_{\nu_i}})
\cong H^*(\lz; 1_{\nu_0}\otimes{\bigotimes}^n_{i=1}
{T^{\lambda_i}_{\nu_i}}^\tconv)\tag1\\
\cong &\,\, H^*(\lz^\tconv; 1_{\nu_0}\otimes
{\bigotimes}^n_{i=1}{T^{\lambda_i}_{\nu_i}}^\tconv) \cong
H^*(\lz^\talg; 1_{\nu_0}\otimes{\bigotimes}^n_{i=1}
{T^{\lambda_i}_{\nu_i}}^\talg).\\
\split
& H^*(\ww; {\bigotimes}^n_{i=1}{T_{\nu_i}})
\cong H^*(\ww; {\bigotimes}^n_{i=1}{T_{\nu_i}}^\tconv)\\
\cong & \,\,H^*(\ww^\tconv; {\bigotimes}^n_{i=1}{T_{\nu_i}}^\tconv)
\cong H^*(\ww^\talg; {\bigotimes}^n_{i=1}{T_{\nu_i}}^\talg).
\endsplit\tag2
\endalign
$$
\endproclaim
\demo{Proof} The second isomorphisms of (1) and (2) follow from [Ka1] 
Theorem 5.3, and the third from ibid\., Corollary A.3.\par
We prove the first of (1) (cf\. ibid\. \S2). 
That of (2) is proved in a similar manner.  
We define the $\lz$ modules $N_\ep$ ($\ep > 0$), $N$ and $N^\tconv$ by 
$$
\split
& N_\ep := H^0(\{\vert z\vert < \ep\}^n; \wnzotln),\\
& N := 1_{\nu_0}\otimes{\bigotimes}^n_{i=1}T^{\lambda_i}_{\nu_i} = 
{\varprojlim}_{\ep\uparrow\infty}N_\ep, \quad \text{and}\\
& N^\tconv := 1_{\nu_0}\otimes{\bigotimes}^n_{i=1}
{T^{\lambda_i}_{\nu_i}}^\tconv = {\varinjlim}_{\ep\downarrow0}N_\ep.
\endsplit
$$
The multiplicative group $\cz$ acts on $N_\ep$  by $T_t: N_\ep \to N_{\ep/\vert t\vert}$ $(t \in \cz)$ as in the proof of Lemma 2.1. Hence $H^*(\lz; N_\ep) = H^*(C^*(\lz; N_\ep)^\ez).$ \par
 By (2.2), if $0 < \delta <1$, the restriction homomorphism 
$C^*(\lz; N_\ep)^\ez \to C^*(\lz; \allowmathbreak N_{\delta\ep})^\ez$ 
has its inverse $T_\delta$, and so is an isomorphism. 
Hence we obtain a series of isomorphisms
$$
\split
& C^*(\lz; N^\tconv)^\ez = \varinjlim C^*(\lz; N_\ep)^\ez\\
= & C^*(\lz; N_\ep)^\ez = \varprojlim C^*(\lz; N_\ep)^\ez 
= C^*(\lz; N)^\ez.
\endsplit
$$
The first isomorphism follows from [K] Lemma 3, p.372 
(see also [Ka1] Lemma 4.3). 
This completes the proof.\qed\enddemo

\demo{Remark 2.4} 
For any $\nu_1 \in \Bbb Z$ the $\ww^\talg$ module 
${T_{\nu_1}}^{\text{formal}} := \cc[[z]]{dz}^{\nu_1}$ of 
$\nu_1$-covariant 
formal tensor fields is just the co-induced module of 
the $\lz^\talg$ module $1_{\nu_1}$. 
This implies ${T_{\nu_1}}^{\text{formal}}\otimes 
\bigotimes^n_{i=2}{T_{\nu_i}}^\talg$ is 
the co-induced module of 
$1_{\nu_1}\otimes\bigotimes^n_{i=2}{T_{\nu_i}}^\talg$.
From the Shapiro lemma (the Frobenius reciprocity)
$$
H^*(\ww^\talg; {T_{\nu_1}}^{\text{formal}}\otimes 
{\bigotimes}^n_{i=2}{T_{\nu_i}}^\talg) 
\cong H^*(\lz^\talg; 1_{\nu_1}\otimes
{\bigotimes}^n_{i=2}{T_{\nu_i}}^\talg).
$$
On the other hand, 
$C^*(\ww^\talg; {T_{\nu_1}}^{\text{formal}}\otimes 
{\bigotimes}^n_{i=2}{T_{\nu_i}}^\talg)^\ez = 
C^*(\ww^\talg; {T_{\nu_1}}^\talg\otimes 
{\bigotimes}^n_{i=2}{T_{\nu_i}}^\talg)^\ez$. 
Hence, by Lemma 2.3, we have
$$
H^*(\ww; {\bigotimes}^n_{i=1}{T_{\nu_i}}) 
\cong H^*(\lz; 1_{\nu_1}\otimes{\bigotimes}^n_{i=2}{T_{\nu_i}}).
$$
This isomorphism can be described in the following explicit way. 
We remark the evaluation map
$$
\ev: T_{\nu_1} \to 1_{\nu_1}, \quad f(z)dz^{\nu_1} \mapsto f(0)1_{\nu_1}
$$
is an $\lz$ homomorphism. So we can define the composite map
$$
H^*(\ww; {\bigotimes}^n_{i=1}{T_{\nu_i}}) 
\to H^*(\lz; {\bigotimes}^n_{i=1}{T_{\nu_i}}) 
\overset\ev\to\to H^*(\lz; 1_{\nu_1}
\otimes{\bigotimes}^n_{i=2}{T_{\nu_i}}),
$$
which coincides with the Shapiro isomorphism.\enddemo

In the succeeding sections we study 
the equivariant cohomology $H^*_\lz(\ccnz; \allowmathbreak 
\wnzotln)$. 
The stalk at $z = (z_1, \dots, z_n) \in \cc^n$ of 
the sheaf $H^*(\lz; \wnzotln)$, 
which appears in the $``E_2$ term (1.1) converging 
to the equivariant cohomology, is given as follows. 
For simplicity we assume, for some $0 \leq r_0 < r_1 < \dots < r_l = n$, 
$$
\cases
z_i = 0, & \quad \text{ if } i \leq r_0\\
z_i = z_{r_k}, & \quad \text{ if } r_{k-1} < i \leq r_k,\\
z_{r_k} \neq 0, & \quad \text{ if }k \geq 1, \quad \text{and}\\
z_{r_k} \neq z_{r_j}, & \quad \text{ if } k \neq j.\\
\endcases\tag2.6
$$
Let $u_k \in H^2(\lz; (T_0)^{\otimes n})$ be the cohomology class 
of a $2$ cocycle defined by
$$
u_k(\xi_1(z)\dz, \xi_2(z)\dz) := 
\sum_{r_{k-1}<i\leq r_k} \int^{z_i}_0 \det\pmatrix
\xi'_1(z) & \xi'_2(z)\\
\xi''_1(z) & \xi''_2(z)\\
\endpmatrix dz
$$
for $\xi_1(z)\dz$ and $\xi_2(z)\dz \in \lz$.
\proclaim{Lemma 2.7}{\rm ([Ka1] Theorem 5.3.)} 
The stalk at $z = (z_1, \dots, z_n) \in \cc^n$ is given by
$$
\multline
H^*(\lz; \wnzotln)_{(z_1, \dots, z_n)}\\
= \cc[u_1, \dots, u_l]\otimes 
H^*(\lz; 1_{\nu_0}\otimes{\bigotimes}^{r_0}_{i=1}T^{\lambda_i}_{\nu_i})
\otimes{\bigotimes}^l_{k=1}H^*(\ww; 
{\bigotimes}_{r_{k-1}<i\leq r_k}T_{\nu_i}).
\endmultline
$$
\endproclaim

\beginsection 3. Fundamental exact sequence.\par
Using the Stein covering $\{U_j\}^n_{j=1}$ of $\ccnz$ given by (1.8), 
we have an $L_0$ isomorphism
$$
H^{n-1}(\ccnz; \wnzotln) = 
1_{\nu_0}\otimes{\bigotimes}^n_{i=1}(T^\times_{\nu_i}/T^{\lambda_i}_{\nu_i}).
$$
From Proposition 1.9 (a corollary of a theorem of Scheja [S]), the $`E_2$ term (1.1) converging to $H^*_\lz(\ccnz; 1_{\nu_0}\otimes\Cal O_{\Bbb C^n}(\allowmathbreak\tau^{\lambda_1,\dots,\lambda_n}_{\nu_1,\dots,\nu_n}))$ is given by 
$$
`E^{p, q}_2 =
\cases
H^p(\lz; 1_{\nu_0}\otimes\bigotimes^n_{i=1}T^{\lambda_i}_{\nu_i}), 
& \quad \text{ if } q = 0,\\
H^p(\lz; \wnztnztln), & \quad \text{ if } q = n - 1,\\
0, & \quad \text{ otherwise. }
\endcases\tag3.1
$$
This means a cohomology exact sequence
$$
\multline
\cdots\to H^{q-n}(\lz; \wnztnztln) \overset{d_n}\to\to 
H^q(\lz; 1_{\nu_0}\otimes\bigotimes^n_{i=1}T^{\lambda_i}_{\nu_i})\\
\to H_\lz^q(\ccnz; \wnzotln)\to\cdots,
\endmultline\tag3.2
$$
which we call {\it the fundamental exact sequence 
for the $\lz$ module 
$1_{\nu_0}\otimes\bigotimes^n_{i=1}T^{\lambda_i}_{\nu_i}$}.
\proclaim{Corollary 3.3} 
If $\nu_0 + \sum^n_{i=1}(\nu_i - \lambda_i) \leq n - 1$, then
$$
H^p(\lz; 1_{\nu_0}\otimes{\bigotimes}^n_{i=1}T^{\lambda_i}_{\nu_i}) 
\cong H_\lz^p(\ccnz; \wnzotln).
$$
\endproclaim
\demo{Proof} The $\lz$ module $\wnztnztln$ is (topologically) 
generated by negative eigenvectors of $\ez =z\dz$ 
under the given assumption. 
Hence $C^*(\lz; \wnztnztln)^\ez = 0$. 
The corollary follows from the sequence (3.2).\qed\enddemo
\par
As an application, we have
\proclaim{Theorem 3.4} 
Let $\nu_0, \nu_1,\dots,\nu_n, \lambda_1,\dots,\lambda_n$ be 
integers satisfying 
$$
\nu_0 + \sum_{i \in I}(\nu_i -\lambda_i) \leq \sharp I - 1
$$
for any non-empty index subset $I \subset \{1, \dots, n\}$. 
Then the restriction homomorphism
$$
H^*(\lz; 1_{\nu_0}\otimes{\bigotimes}^n_{i=1}T^{\lambda_i}_{\nu_i}) 
\to H^*(\lz; 1_{\nu_0}\otimes{\bigotimes}^n_{i=1}T^\times_{\nu_i})
$$
is an isomorphism.
\endproclaim
\demo{Proof} 
We prove the theorem by induction on $n$. 
We have $C^*(\lz; 1_{\nu_0}\otimes T^{\lambda_1}_{\nu_1})^\ez 
= C^*(\lz; 1_{\nu_0}\otimes T^\times_{\nu_1})^\ez$, 
which implies the theorem for $n = 1$.\par
If $n \geq 2$, then, by Corollary 3.3, we have 
$$
H^*(\lz; 1_{\nu_0}\otimes{\bigotimes}^n_{i=1}T^{\lambda_i}_{\nu_i}) 
\cong H^*_\lz(\ccnz; \wnzotln).
$$
Consider the inclusion map $i: (\cz)^n \to \ccnz$, 
where $\cz = \cc - \{0\}$.  
We calculate the stalk at the point $z \in (\ccnz) - (\cz)^n$ 
of the sheaf $H^*(\lz; i_*\wnzoztn)$. 
For simplicity we assume the point $z$ satisfies the condition (2.6). 
Lemma 2.7 (or [Ka1] Theorem 5.3) implies
$$
\multline
H^*(\lz; i_*\wnzoztn)_{(z_1, \dots, z_n)}\\
= \cc[u_1, \dots, u_l]\otimes 
\varinjlim_{\ep\downarrow 0}(H^*(\lz; 1_{\nu_0}
\otimes\Cal O_{(\cz)^n}(\tau_{\nu_1,\dots,\nu_{r_0}})
(\{0<\vert z\vert<\ep\}^{r_0})))\\
\otimes{\bigotimes}^l_{k=1}
H^*(\ww; {\bigotimes}_{r_{k-1}<i\leq r_k}T_{\nu_i}).
\endmultline
$$
The second term of the RHS is isomorphic to 
$H^*(\lz; 1_{\nu_0}\otimes\bigotimes^{r_0}_{i=1}T^\times_{\nu_i})$ 
by a similar method to the proof of Lemma 2.3. 
We have $r_0 < n$ since $z \in \ccnz$.
By the inductive assumption the term is isomorphic 
to $H^*(\lz; 1_{\nu_0}\otimes\bigotimes^{r_0}_{i=1}T^{\lambda_i}_{\nu_i})$. 
Hence we have an isomorphism of sheaves over $\ccnz$
$$
H^*(\lz; i_*\wnzoztn) \cong H^*(\lz; \wnzotln)\vert_\ccnz.
$$
Comparing the $``E_2$ terms, we obtain
$$
H^*_\lz(\ccnz; \wnzotln)\cong H^*_\lz(\ccnz; i_*\wnzoztn).
$$
We write simply $\Cal F := \wnzoztn$. 
For any open subset $U$ of $\ccnz$ we have 
$$
C^*(\lz; i_*\Cal F)(U) = C^*(\lz; (i_*\Cal F)(U)) = 
C^*(\lz; \Cal F(U\cap (\cz)^n)) = i_*C^*(\lz; \Cal F)(U).
$$ 
Hence $C^*(\lz; i_*\Cal F) = i_*C^*(\lz; \Cal F)$. 
If $U$ is Stein, then $U \cap (\cz)^n$ is also Stein, 
so that $H^q(U \cap (\cz)^n; C^*(\lz; \Cal F)) = 0$ 
for $q \neq 0$ from Proposition 1.5. 
This implies $\Cal R^qi_*C^*(\lz; \Cal F) = 0$ for $q \neq 0$, 
so that 
$H^*(\ccnz; i_*C^*(\lz; \Cal F)) = H^*((\cz)^n; \allowmathbreak 
C^*(\lz; \Cal F))$.
Comparing the $`E_2$-terms, we have 
$$
H^*_\lz(\ccnz; i_*\Cal F) =H^*_\lz((\cz)^n; \Cal F) 
= H^*(\lz; 1_{\nu_0}\otimes{\bigotimes}^n_{i=1}T^\times_{\nu_i}).
$$
This completes the induction.\qed\enddemo
We introduce an algebraic variant of the $\lz$ module 
${T^\times_\nu}$ by
$$
{T^\times_\nu}^\talg := \varinjlim_{\lambda\uparrow\infty}
{T^\lambda_\nu}^\talg = 
\varinjlim_{\lambda\uparrow\infty}\frac1{z^\lambda}\cc[z]dz^\nu 
= \Bbb C[z^{-1}, z]dz^\nu.
$$
The $\lz^\talg$ module $1_{\nu_0}\otimes\bigotimes^n_{i=1}
{T^\times_{\nu_i}}^\talg$ is the space of algebraic tensors 
on the complex torus $(\cc - \{0\})^n$. 
It follows from Theorem 3.4 and Lemma 2.3

\proclaim{Corollary 3.5}
$$
H^*(\lz; 1_{\nu_0}\otimes{\bigotimes}^n_{i=1}T^\times_{\nu_i}) 
= H^*(\lz^\talg; 1_{\nu_0}\otimes{\bigotimes}^n_{i=1}
{T^\times_{\nu_i}}^\talg).
$$
\endproclaim
\par
We denote by $F(\Omega)$ the Fr\'echet space of complex analytic 
functions on a complex analytic manifold $\Omega$ with the topology 
of uniform convergence on compact sets.  
For any open subset $U \subset \cc^n$ the Lie algebra $\ww$ 
acts on the Fr\'echet space $F(U)$ by the diagonal Lie derivative action. 
Clearly we have an $\lz$ isomorphism $F(\cc^n) \cong (T_0)^{\otimes n}$. 
Now we determine the cohomology algebra $H^*(\ww; F(\cc^n))$ explicitly 
by making a comparison with that with values in the functions 
on the configuration space
$$
P_n = \{(z_1, \dots, z_n) \in \Bbb C^n; z_i \not= z_j 
\,\, (i \not=j)\} \quad(n \ge 1).
$$

\proclaim{Theorem 3.6} The restriction homomorphism
$$
\iota_n: H^*(\ww; F(\cc^n)) \to H^*(\ww; F(P_n))
$$
is isomorphic for any $n \ge 1$.
\endproclaim
It follows from [Ka1] Proposition 9.8
\proclaim{Corollary 3.7}
$$
H^*(\ww; F(\cc^n)) = 
H^*(P_n)\otimes\cc[v_2, \dots, v_n]
\otimes{\bigwedge}^*(\nabla^1_0, \dots, \nabla^n_0).
$$
Here the (holomorphic) de Rham cohomology algebra of $P_n$, 
$H^*(P_n)$, is mapped into the algebra $H^*(\ww; F(\cc^n))$ 
in an obvious way (ibid., Corollary 9.9). The $2$ cocycle $v_i$ is defined by 
$$
v_i(\xi_1(z)\dz, \xi_2(z)\dz) := \int^{z_i}_{z_{i-1}}\det\pmatrix
\xi'_1(z) & \xi'_2(z)\\
\xi''_1(z) & \xi''_2(z)\\
\endpmatrix dz
$$
and the $1$ cocycle $\nabla^j_0$ by
$$
\nabla^j_0(\xi(z)\dz) = \xi'(z_j)
$$
for $\xi(z)\dz, \xi_1(z)\dz \text{ and } \xi_2(z)\dz \in \lz$.
\endproclaim
\demo{Proof of Theorem 3.6} 
We prove the theorem by induction on $n$. 
It is already known for the case $n \le 2$ [FF]. 
From [Ka1], Corollary 9.9, the homomorphism $\iota_n$ is 
surjective for any $n \geq 1$. \par
Suppose $n \ge 2$. 
We denote $P^\times_n := P_n \cap (\cc - \{0\})^n$. 
Then we have a commutative diagram
$$
\CD
H^*(\ww; F(\cc^{n+1})) @>{\iota_{n+1}}>> H^*(\ww; F(P_{n+1})\\
@VVV @VVV\\  
H^*(\lz; F(\cc^{n})) @>{\iota'}>> H^*(\lz; F(P^\times_n)),
\endCD
$$
where the vertical homomorphisms are the Shapiro homomorphisms 
given in Remark 2.4, and the lower horizontal $\iota'$ 
the restriction homomorphism.
Since the left vertical is isomorphic from Remark 2.4 and 
$\iota_{n+1}$ is surjective, it suffices to show $\iota'$ 
is isomorphic. 
From Corollary 3.3  follows 
$H^*_\lz(\cc^n - \{0\}; \Cal O_{\cc^n}) = H^*(\lz; F(\cc^n)).$\par
Consider the inclusion map $i: P^\times_n \to \ccnz$. 
If $(z_1,\dots, z_n) \in \ccnz$, 
all the multiplicities of the set $\{ 0, z_1, \dots, z_n\}$ 
are not greater than $n$. 
Hence we deduce an isomorphism of sheaves over $\ccnz$
$$
H^*(\lz; i_*\Cal O_{P^\times_n}) 
\cong H^*(\lz; \Cal O_{\cc^n})\vert_{\ccnz}
$$
from the inductive assumption together with a similar method 
to the proof of Theorem 3.4. 
Moreover we have
$$
H^*_\lz(\cc^n - \{0\}; i_*\Cal O_{P^\times_n}) 
\cong H^*_\lz(P^\times_n; \Cal O_{P^\times_n}) 
= H^*(\lz; F(P^\times_n))
$$
because $P^\times_n$ is Stein. 
Consequently the restriction homomorphism 
$\iota'$ is an isomorphism, which completes the induction.\enddemo

\beginsection 4. Action of the symmetric group.\par

The $n$-th symmetric group $\sn$ acts on the $\lz$ module $\wnztlnn$ 
in two ways:  
If $\sigma \in \sn$ and $\wnzo f(z_1,\dots, z_n)\dzn \in \wnztlnn$, 
then $\sigma(\wnzo f(z_1,\dots, z_n)\dzn)$ is given by
$$
\wnzo f(z_{\sigma(1)},\dots, z_{\sigma(n)})\dzn \in \wnztlnn
$$
in one way, and by
$$
\sgns\wnzo f(z_{\sigma(1)},\dots, z_{\sigma(n)})\dzn \in \wnztlnn
$$
in the other way.  
We call the former {\it the symmetric action} and 
the latter {\it the alternating action} of the group $\sn$.  
The $\lz$ module $\wnztnztlnn$ has the symmetric and 
the alternating actions of the group $\sn$ in a similar manner.\par
\proclaim{Proposition 4.1} 
The fundamental exact sequence (3.2) for the $\lz$ module $\wnztlnn$ 
is $\sn$ equivariant, 
if the group $\sn$ acts on $H^*(\lz; \wnztnztlnn)$ 
by the symmetric (resp\. alternating) action and 
on the other two terms by the alternating 
(resp\. symmetric) action.
\endproclaim
\demo{Proof} 
Consider the alternating \v Cech complex 
$C^*(\frak U) = C^*(\frak U; \wnzotlnn)$ 
with respect to the Stein covering $\frak U = \{U_j\}^n_{j=1}$ 
of $\ccnz$ given by (1.8).  
The cohomology group of the total complex of the double complex
$$
C^*(\lz; C^*(\frak U; \wnzotlnn))
$$
is the equivariant cohomology group $H^*_\lz(\ccnz; \wnzotlnn)$.\par
The symmetric group $\sn$ acts on the \v Cech complex as follows.  
Suppose a $q$ cochain $c \in C^q(\frak U)$ is given by the tensors
$$
\wnzo c_{i_0\dots i_q}(z_1,\dots, z_n)\dzn 
\quad\quad(1 \leq i_0, \dots, i_q \leq n),
$$
where $c_{i_0\dots i_q}(z_1,\dots, z_n)\dzn \in 
\Cal O_{\cc^n}((\tau^\lambda_\nu)^n)(\bigcap^q_{\alpha=0}U_{i_\alpha})$. 
If $\sigma \in \sn$, the cochain $\sigma c$ is given by the tensors
$$
\split
& \wnzo (\sigma c)_{i_0\dots i_q}(z_1,\dots, z_n)\dzn,\\
& (\sigma c)_{i_0\dots i_q}(z_1,\dots, z_n) := 
\sgns c_{\sigma^{-1}(i_0)\dots \sigma^{-1}(i_q)}
(z_{\sigma(1)},\dots, z_{\sigma(n)}).
\endsplit
$$
Then the group $\sn$ acts on $H^0(C^*(\frak U)) = \wnztlnn$ 
by the alternating action and on $H^{n-1}(C^*(\frak U))= \wnztnztlnn$ 
by the symmetric action, as was to be shown.  
The rest is proved in a similar way.\qed\enddemo
\par

To put the proposition into practice, we review on sheaves 
on which a finite transformation group acts. 
Here we regard a sheaf as an etale covering space.\par
Let $M$ be a Hausdorff space on which a finite group $G$ 
acts continuously. 
The isotropy group at a point $x \in M$ is denoted by $G_x$, 
i.e., $G_x := \{ \gamma \in G ; \,\,\gamma x = x \}$. 
 Let $\pi: \fai \to M$ be a sheaf of complex vector spaces 
on which the group $G$ acts continuously and linearly. 
The subset $\faig$ of $\fai$ defined by 
$\faig := \{a \in \fai; \forall\gamma \in G_{\pi(a)},\,\, \gamma a = a\}$ 
is a linear subsheaf of $\fai$. 
The functor $\fai \mapsto \faig$ is exact. 
Furthermore we have a canonical isomorphism
$$
\Cal C^*(M; \faig) \overset\cong\to\to \Cal C^*(M; \fai)^G,\tag4.2
$$
where $\Cal C^*(M; \cdot)$ denotes the canonical resolution. 
See, for example, [B]II\S2. The isomorphism (4.2) follows from the fact
$$
\faig(U) = \{s \in \fai(U); \,\,\forall x\in U, \forall\gamma \in G_x,
\,\, \gamma s(x) = s(x)\} \tag 4.3
$$
and the definition of the resolutiom $\Cal C^*(M; \cdot)$.\par
If $\varpi: M \to M/G$ is the quotient map, we have
$$
\split
& H^*(M; \fai)^G = H^*(M/G; \varpi_*\fai)^G 
= H^*(\Cal C^*(M/G; \varpi_*\fai)(M/G)^G)\\
= & H^*(\Cal C^*(M/G; \varpi_*\fai)^G(M/G)) 
= H^*(\Cal C^*(M/G; (\varpi_*\fai)^G)(M/G))\\
= & H^*(M/G; (\varpi_*\fai)^G).
\endsplit
$$
The third equality follows from (4.3) and the fourth from (4,2). 
The quotient space $\faig/G$ is a sheaf over $M/G$ and satisfies
$$
(\faig/G)(U) = \fai(\varpi^{-1}(U))^G = (\varpi_*\fai)^G(U)
$$
for any open subset $U$ of $M/G$. 
This means $\faig/G = (\varpi_*\fai)^G$. 
Consequently we obtain an isomorphism
$$
H^*(M; \fai)^G \cong H^*(M/G; \faig/G).\tag4.4
$$
\par
For a topological vector space $V$, 
the completed $n$ fold symmetric (resp\. alternating) tensor product 
is denoted by $S^n(V)$ (resp\. $\bigwedge^nV$).  
We denote by $(\cdot)^{\sn,\tsym}$ (resp\. $(\cdot)^{\sn,\talt}$) 
the space of invariants of the space $\cdot$ under the symmetric 
(resp\. alternating) action of the group $\sn$.  
We have $S^n(V) = (V^{\otimes n})^{\sn,\tsym}$ and 
$\bigwedge^nV = (V^{\otimes n})^{\sn,\talt}$. \par
Proposition 4.1 implies a cohomology exact sequence
$$
\multline
\cdots\to H^{q-n}(\lz; 1_{\nu_0}\otimes S^n(T^\times_\nu/T^\lambda_\nu)) 
\overset{d_n}\to\to H^q(\lz; 1_{\nu_0}\otimes{\bigwedge}^nT^\lambda_\nu)\\
\to H_\lz^q(\ccnz; 
1_{\nu_0}\otimes\Cal O_{\cc^n}((\tau^\lambda_\nu)^n))^{\sn,\talt}
\to\cdots,
\endmultline\tag4.5
$$
which we call {\it the fundamental exact sequence for 
the $\lz$ module $1_{\nu_0}\otimes\bigwedge^nT^\lambda_\nu$}. 
From (4.4) the $``E_2$ term (1.1) converging 
to the third term of (4.5) is given by
$$
``E^{p, q}_2 = H^p(\ccnz/\sn; H^q(\lz; \wnzotlnn)^{\sn,\talt}/\sn).
$$
By Lemma 2.7 we have
\proclaim{Lemma 4.6} If $z = (z_1,\dots,z_n) \in \cc^n$ satisfies 
the condition (2.6), 
the stalk at $z$ of the sheaf $H^*(\lz; \wnzotlnn)^{\sn,\talt}$ 
is given by 
$$
\multline
H^*(\lz; \wnzotlnn)^{\sn,\talt}_{(z_1, \dots, z_n)}\\
= \cc[u_1, \dots, u_l]\otimes 
H^*(\lz; 1_{\nu_0}\otimes{\bigwedge}^{r_0}T^{\lambda}_{\nu})\otimes
{\bigotimes}^l_{k=1}H^*(\ww; {\bigwedge}^{r_k - r_{k-1}}T_\nu).
\endmultline
$$
\endproclaim
Similar results to (4.5) and Lemma 4.6 hold 
for the $\lz$ module $1_{\nu_0}\otimes S^nT^\lambda_\nu$.\par
\medskip

Finally we prove a certain vanishing theorem for the cohomology of $\ww$ 
with values in the quadratic differentials on $\cc^n$.  
For the rest we denote $Q = T_2$, $Q^1 = T^1_2$ 
and $Q^\times = T^\times_2$.  
The cohomology groups $H^q(\ww; \bigwedge^pQ)$ and 
$H^q(\lz; \bigwedge^pQ^1)$ are closely related 
to the $(p, q)$ cohomology groups of the parameter (moduli) space 
and the total space of the universal family of dressed Riemann surfaces 
respectively [ADKP] [Ka2].
\proclaim{Theorem 4.7}
$$
\split
& H^q(\lz; 1_2\otimes{\bigwedge}^{n-1}Q) = 0,  \quad\text{ if } q < n,\\
& H^q(\ww; {\bigwedge}^nQ) = 0,  \quad\text{ if } q < n.
\endsplit
$$
\endproclaim
\demo{Proof} We remark the first part implies the second part for each $n$. 
The first part is already proved in [FF] for $n \leq 2$.\par
Suppose $ n \geq 2$.  If $q \leq n$, we have 
$$
H^q(\lz; 1_2\otimes\Cal O_{\cc^n}((\tau_2)^n))^{\sn, \talt}\vert_{\ccnz} = 0
$$
by the inductive assumption and Lemma 4.6.  
It follows from the fundamental exact sequence (4.5)
$$
H^q(\lz; 1_2\otimes{\bigwedge}^nQ) = 
\cases
H^0(\lz; 1_2\otimes S^n(Q^\times/Q)), & \text{ if $q = n$,}\\
0, & \text{ if $q < n$.}
\endcases
$$
Hence the following lemma completes the induction. 
Here we denote by $q_\nu$ the class $(z^{\nu-2}dz^2\mod Q) \in Q^\times/Q$. 
Clearly $\ez q_\nu = \nu q_\nu$.
\proclaim{Lemma 4.8} 
$$
H^0(L_1; S^*(Q^\times/Q)) = \cc[q_1] 
\quad\text{(the polynomial algebra in $q_1$)},
$$
where $L_1$ denotes the subalgebra of $\lz$ generated by 
$e_k := z^{k+1}\dz$ ($k \geq 1$).\endproclaim
\demo{Proof of Lemma 4.8} 
Assume there exists some $\alpha \in S^m(Q^\times/Q)^{L_1} - \cc {q_1}^m$, 
$m \geq 1$. 
Let $q_{-p}$, $p \geq 0$, be the maximal 
among those which appear in $\alpha$ except $q_1$. 
Then we have
$$
\alpha = {\sum}_{i+j \leq m} {q_1}^i{q_{-p}}^j f_{ij}(q_{-p-1}, \dots, 
q_{-s}, \dots).
$$
There exists $i_0 := \max\{i; \exists j \geq 1, \,\, 
f_{ij}\neq 0\}$ because of the choice of $q_{-p}$. 
Now we have 
$$
\split
0 = \,e_{p+1}\cdot\alpha
= & \,{\sum}_{i+j \leq m} j(p+2){q_1}^{i+1}{q_{-p}}^{j-1}
 f_{ij}(q_{-p-1}, \dots, q_{-s}, \dots)\\
& + {\sum}_{i+j \leq m} {q_1}^i{q_{-p}}^j e_{p+1}( 
f_{ij}(q_{-p-1}, \dots, q_{-s}, \dots)).
\endsplit
$$
Since $e_{p+1}(f_{ij})$ has no $q_1$'s, 
we have 
$$
{\sum}_{j\leq m-i_0}j{q_1}^{i_0+1}{q_{-p}}^{j-1}f_{i_0j}(q_{-p-1}, \dots, 
q_{-s}, \dots) = 0,
$$
or equivalently, $f_{i_0j} = 0$ for any $j \geq 1$. 
This contradicts the choice of $i_0$. 
Hence $S^*(Q^\times/Q)^{L_1} = \cc[q_1]$, as was to be shown.
\qed
\enddemo
Consequently the proof of Theorem 4.7 is completed.
\qed\enddemo
\proclaim{Theorem 4.9}
$$
H^q(\lz; {\bigwedge}^nQ^1) = 0, \quad \text{if } q < n.
$$
\endproclaim
\demo{Proof} The theorem for $n = 1$ is shown in [FF].\par
Suppose $n \geq 2$. By the inductive assumption and Theorem 4.7 
we have
$$
H^q(\lz; \Cal O_{\cc^n}((\tau^1_2)^n))^{\sn, \talt}\vert_\ccnz = 0
$$
for $q < n$.  Hence, by the fundamental exact sequence (4.5), 
$H^q(\lz; \bigwedge^nQ^1) = 0$ for $q < n$, as was to be shown.\qed\enddemo
\demo{Remark 4.10} \roster
\item We have $H^n(\ww; \bigwedge^nQ) \neq 0$ and 
$H^n(\lz; \bigwedge^nQ^1) \neq 0$.  
In fact, if $\nabla_2$ is the $1$ cocycle in $C^1(\ww; Q)$ defined by 
$$
\nabla_2(\xi(z)\dz) = \xi^{(3)}(z)dz^2, \quad (\xi(z)\dz \in \ww),
$$
then the cohomology class of the $n$-th power $(\nabla_2)^n$ 
does not vanish in $H^n(\ww; \bigwedge^nQ)$ nor $H^n(\lz; \bigwedge^nQ^1)$.  
Hence the estimates in Theorems 4.7 and 4.9 are the best possibles.
\item In Theorem 4.7 it is essential to take the {\it alternating} 
tensor product. 
By a computation due to Kazushi Ahara 
we have $\dim H^4(\ww; Q^{\otimes5}) = 3$.
\item Theorems 4.7 and 4.9 suggest that the Virasoro equivariant $(p, q)$ 
cohomology of the parameter space and the total space of 
the universal family of dressed Riemann surfaces 
would vanish for $p > q$ [Ka2].
\endroster
\enddemo

\par
\vskip 15mm
\noindent
Department of Mathematical Sciences,
\newline
University of Tokyo
\newline
Tokyo, 153-8914 Japan 
\newline
e-mail address: kawazumi\@ms.u-tokyo.ac.jp
\par

\end